\documentclass[12pt]{amsart}
\usepackage{latexsym, amssymb, amsmath, amscd}
 \usepackage{eucal}

    \newtheorem{rema}{Remark}[section]
    \newtheorem{propo}[rema]{Proposition}
\newtheorem{algo}[rema]{Algorithm}
   \newtheorem{theo}[rema]{Theorem}
   \newtheorem{def-theo}[rema]{Definition-Theorem}

    \newtheorem{lemma}[rema]{Lemma}
    \newtheorem{corol}[rema]{Corollary}
     \newtheorem{exam}[rema]{Example}
  \newtheorem{rmk}[rema]{Remark}

	\newcommand{\nno}{\nonumber}

	\newcommand{\p}{\partial}

 \newcommand{\res}{\mbox{\rm Res}}

 \newcommand{\pf}{{\it Proof:}\hspace{2ex}}
 
 \newcommand{\epfv}{\hspace{1em}$\Box$\vspace{1em}}

\newcommand{\bC}{{\mathbb C}}
\newcommand{\bZ}{{\mathbb Z}}

\newcommand{\bN}{{\mathbb N}}
\newcommand{\bT}{{\mathbb T}}

\newcommand{\cP}{{\mathcal P}}

\newcommand{\cDzz}{{\mathcal D \langle \langle z \rangle\rangle}}
\newcommand{\cDkzz}{{\mathcal D_k \langle \langle z \rangle\rangle}}

\newcommand{\cDrzz}{{\mathcal D er\langle \langle z \rangle\rangle}}

\newcommand{\cDrkzz}{{\mathcal D er_k \langle \langle z \rangle\rangle}}

\newcommand{\cNcs}{{${\mathcal N}$CS} }

\newcommand{\kz}{{k\langle z \rangle}}
\newcommand{\kzz}{{k\langle \langle z \rangle\rangle}}
\newcommand{\ktz}{{k[t]\langle z \rangle}}
\newcommand{\ktzz}{{k[t]\langle \langle z \rangle\rangle}}
\newcommand{\kttzz}{{k[[t]]\langle \langle z \rangle\rangle}}

\newcommand{\BQ}{\begin{eqnarray}}
\newcommand{\EQ}{\end{eqnarray}}
\newcommand{\BQn}{\begin{eqnarray*}}
\newcommand{\EQn}{\end{eqnarray*}}
\newcommand{\BL}{\begin{align}}
\newcommand{\EL}{\end{align}}
\newcommand{\BLn}{\begin{align*}}
\newcommand{\ELn}{\end{align*}}
\newcommand{\BA}{\begin{align}}
\newcommand{\EA}{\end{align}}
\newcommand{\BAn}{\begin{align*}}
\newcommand{\EAn}{\end{align*}}
\newcommand{\wtilde}{\widetilde}

\newcommand{\lb}{\left[}
\newcommand{\rb}{\right]}
\newcommand{\lp}{\left(}
\newcommand{\rp}{\right)}

\newcommand{\fr}{\frac}
\newcommand{\pz}{\frac{\p}{\p z}}

\title[Deformations and Inversion Formulas] 
{Deformations and Inversion Formulas for 
Formal Automorphisms in Noncommutative Variables}

    \author{Wenhua Zhao}
      
 \thanks{{\it Address}: {Department of Mathematics, Illinois State 
 University, Normal, IL 61790-4520.}  
\quad {\it E-mail}: {wzhao@ilstu.edu} }
\date{\today}

\begin{document}

\begin{abstract}
Let $z=(z_1, z_2, ... , z_n)$ be noncommutative 
free variables and $t$ a formal parameter
which commutes with $z$. Let $k$ be any unital 
integral domain of any characteristic and   
$F_t(z)=z-H_t(z)$ with 
$H_t(z)\in {k[[t]]\langle \langle z \rangle\rangle}^{\times n}$ 
and the order $o(H_t(z))\geq 2$.
Note that $F_t(z)$ can be viewed 
as a deformation of the formal map 
$F(z):=z-H_{t=1}(z)$ when it makes sense 
(for example, when 
$H_t(z)\in {k[t]\langle \langle z \rangle\rangle}^{\times n}$).  
The inverse map 
$G_t(z)$ of $F_t(z)$ can always be written 
as $G_t(z)=z+M_t(z)$
with $M_t(z)\in {k[[t]]\langle \langle z \rangle\rangle}^{\times n}$ 
and $o(M_t(z))\geq 2$. 
In this paper, we first derive the PDE's 
satisfied by $M_t(z)$ and 
$u(F_t), u(G_t)\in {k[[t]]\langle \langle z \rangle\rangle}$ 
with $u(z)\in {k\langle \langle z \rangle\rangle}$ 
in the general case 
as well as in the special case when $H_t(z)=tH(z)$ 
for some $H(z)\in {k\langle \langle z \rangle\rangle}^{\times n}$. 
We also show that the formal power series above are
actually characterized by 
certain Cauchy problems
of these PDE's. Secondly,
we apply the derived PDE's 
to prove a recurrent inversion formula 
for formal maps in noncommutative variables. 
Finally, for the case char.\,$k=0$, we derive
an expansion inversion formula by 
the planar binary rooted trees.   


\end{abstract}


 \bibliographystyle{alpha}
    \maketitle
 
\renewcommand{\theequation}{\thesection.\arabic{equation}}
\renewcommand{\therema}{\thesection.\arabic{rema}}
\setcounter{equation}{0}
\setcounter{rema}{0}
\setcounter{section}{0}

\section{\bf Introduction} \label{S1}

Let $z=(z_1, z_2, ... , z_n)$ 
be $n$ noncommutative free 
variables and $t$ a formal parameter
which commutes with $z$. We fix a unital 
integral domain $k$ of any characteristic
and denote by 
$\kzz$ and $\kttzz$ the algebras 
of formal power series in 
$z$ over $k$ and $k[[t]]$, 
respectively. 
In this paper, we first study 
the deformations of automorphisms 
of $\kzz$ parameterized by 
the formal parameter $t$ and 
then derive some inversion formulas 
for the automorphisms of $\kzz$.
More precisely, we consider the 
automorphisms $F_t(z)$ of $\kttzz$ 
over $k[[t]]$ of the form 
$F_t(z)=z-H_t(z)$ with 
$H_t(z)\in \kttzz^{\times n}$ 
and the order $o(H_t(z))\geq 2$. 
Note that
$F_t(z)$ can be viewed 
as a general deformation 
parameterized by $t$ of 
the formal map 
$F(z):=z-H_{t=1}(z)$ 
when it exists (for example, 
when $H_t(z)\in \ktzz^{\times n}$, or 
when $k=\bC$ and all coefficients of $H_t(z)$ 
are holomorphic functions of $t$ which are 
convergent over an open subset of $\bC$ 
containing the closed unit disk).
In particular, this is indeed 
the case for the special 
deformation $F_t(z)=z-tH(z)$ with 
$H(z)\in \kzz^{\times n}$, 
i.e. $H_t(z)=tH(z)$. 
We will always denote by $G_t(z)$ the formal 
inverse map of $F_t(z)$ and write it 
as $G_t(z)=z+M_t(z)$
with $M_t(z)\in \kttzz^{\times n}$ 
and $o(M_t(z))\geq 2$.
When $F_t(z)$ is 
the special deformation 
$F_t(z)=z-tH(z)$ above, 
we also write its inverse as 
$G_t(z)=z+tN_t(z)$ with 
$N_t(z)\in \kttzz^{\times n}$. 
In the first part of this paper, 
we derive the PDE's in $z$ and $t$  
satisfied by $M_t(z)$, $N_t(z)$ $u(F_t)$ and $u(G_t)$ 
$(u(z)\in \kzz)$. 
In particular, 
we show that $N_t(z)$ 
is a formal power series
of the Cauchy problem of a
Burgers-like PDE (see Theorem \ref{T4.3} 
and Remark \ref{rmk-Burgers}).
When $char.\, k=0$,
$N_t(z)$ is actually the unique power 
series solution 
of a Cauchy problem of the PDE; 
while when $char.\, k=p>0$,  
$N_t(z)$ is completely 
determined by this property 
together with its coefficients of $t^{mp}$ 
$(m\geq 1)$, which can be calculated 
by some other methods 
(see Corollary \ref{Uniqueness} 
and Theorem \ref{T5.5}). In addition, 
we also discuss some other characterizing 
properties of $N_t(z)$.
In the second part of this paper,
we apply the PDE satisfied by 
$N_t(z)$ to 
derive a recurrent 
inversion formula 
and, when $char.\, k=0$, 
an expansion inversion formula 
by the planar binary rooted trees 
for formal maps in 
noncommutative 
free variables.
Note that the special 
deformation 
$F_t(z)=z-tH(z)$ 
for commutative variables $z$
over any unital commutative 
ring $k$ of characteristic 
zero has been studied 
in \cite{BurgersEq}. 
Here we not only generalize 
the results in 
\cite{BurgersEq} to formal 
maps in noncommutative 
variables, but also give
some inversion algorithms
for the case when the base ring $k$ 
has $char.\, k=p>0$. 
When $char.\, k=0$, 
the expansion 
inversion formula by the planar 
binary rooted trees 
for the symmetric maps 
in \cite{BurgersEq} 
is also generalized to 
general automorphisms.

The problem seeking various 
inversion formulas of formal maps in 
commutative variables has a long 
history in mathematics. 
It started with the 
Lagrange's inversion  
formula in one variable by L. Lagrange \cite{L} in 1770, 
then the Jacobi's inversion formula by C. G. J. Jacobi
\cite{J1} in $1830$ and \cite{J2} in $1844$. Later, 
motivated by the well-known Jacobian conjecture 
proposed by O. H. Keller \cite{Ke} in $1939$,
more inversion formulas 
have been proved (see \cite{BCW}, 
\cite{E}, \cite{Sm} and references there 
for more history and known results 
on the Jacobian conjecture).
In $1965$, I. G. Good \cite{Go} generalized 
the Lagrange's inversion formula 
to the multiple variable case.
In $1974$,  Gurjar (unpublished) 
and later Abhyankar \cite{Ab} proved 
so-called Abhyankar-Gurjar inversion formula. 
In $1981$, 
H. Bass, E. Connell and D. Wright \cite{BCW} 
and D. Wright \cite{Wr1} 
proved the so-called Bass-Connell-Wright's 
tree expansion formula. 
Very recently, D. Wright and the author \cite{WZ}
generalized this formula
to tree expansion formulas for the D-log and 
the formal flow of formal maps. 
In \cite{BurgersEq} and \cite{HNP}, 
the author proved 
a recurrent inversion formula 
in general and a non-recurrent 
formula for the symmetric maps 
which satisfy the Jacobian condition. 
The later was mainly 
motivated by the remarkable symmetric 
reduction on the Jacobian conjecture 
achieved recently by  M. de Bonlt 
and A. van den Essen in \cite{BE1} 
and G. Meng in \cite{M}. 

On the other hand, comparing with the 
commutative case, 
it seems not many inversion formulas for 
formal automorphisms in noncommutative  
variables are known in the literature. 
But, for an interesting approach to 
this problem, see \cite{Ge1};  
for several $q$-analogue inversion formulas  
see \cite{An}, \cite{Ga}, \cite{GH}. 

One remark is that, based on some results 
obtained in this paper, later, in the 
followed papers \cite{GTS-II}, \cite{GTS-III} and 
\cite{GTS-V}, some connections of 
the commutative or noncommutative inversion
problem with the Hopf algebra ${\mathcal N}Sym$ 
of noncommutative symmetric functions, which were 
first introduced and studied in \cite{G-T},  
and the Grossman-Larson Hopf algebra (\cite{GL}, \cite{F}) 
of labeled rooted trees will be studied. 
In particular, more inversion formulas 
in both commutative and 
noncommutative cases will 
be derived in \cite{GTS-III}. 
The tree expansion formulas 
obtained in \cite{BCW}, \cite{Wr1} 
and \cite{WZ} for the inverse map, 
the D-Log's and the formal flows 
in the commutative case will 
also be generalized 
in \cite{GTS-V} to 
the noncommutative case.

The arrangement of this paper is as follows. 
In Section \ref{S2}, 
we first fix some notation which will be 
used throughout the paper. We then 
consider certain properties 
of derivations and differential operators 
in noncommutative variables.
In particular, we prove two chain rules 
for the derivations of $\kzz$ and $\kttzz$, 
respectively (see Lemma \ref{k-ChainRule} and \ref{t-ChainRule}).
In Section \ref{S3} and \ref{S4}, we study 
the general deformation 
$F_t(z)=z-H_t(z)$ with $H_t(z)\in \kttzz^{\times n}$
and the special deformation
$F_t(z)=z-tH(z)$ 
with $H(z) 
\in \kzz^{\times n}$, respectively.
We not only derive the PDE's satisfied by 
$M_t(z)$, $N_t(z)$ as well as 
formal power series of 
the forms $u(F_t)$ and $u(G_t)$ 
with $u(z)\in \kzz$, 
but also show that
the elements above are 
also characterized by 
certain Cauchy problems 
of these PDE's.
Note that, not all these results 
are needed later for the derivations 
of the inversion formulas 
in the second part of this paper, 
but will be crucial for 
the followed papers \cite{GTS-II}, 
\cite{GTS-III} and \cite{GTS-V}.  
In Section \ref{S5}, 
we apply some results 
obtained in Section \ref{S4} to derive 
a recurrent inversion formula 
for formal maps in noncommutative 
variables over a base ring $k$ 
of any characteristic. 
In Section \ref{S6}, we assume 
our base ring $k$ has characteristic zero 
and prove an expansion inversion 
formula by the planar binary 
rooted trees. 

One final remark is as follows. For simplicity, 
we mainly focus on formal maps in 
noncommutative variables $z$. 
But most of the results obtained in this paper
have their analogs for commutative variables, 
which can be derived either 
by taking the quotient over 
the ideal generated by the commutators of 
$z_i$'s or by applying parallel arguments.

{\bf Acknowledgment}:
The author would like to thank the referee 
for pointing out misprints and providing 
valuable suggestions.

\renewcommand{\theequation}{\thesection.\arabic{equation}}
\renewcommand{\therema}{\thesection.\arabic{rema}}
\setcounter{equation}{0}
\setcounter{rema}{0}

\section{\bf Chain Rules in the Noncommutative Case} \label{S2}

In this section,  
we consider 
certain properties 
of derivations and 
differential operators 
in noncommutative free variables.
In particular, we prove two variations 
of the usual chain rule in the commutative case
for the derivations in the noncommutative case
(see Lemma \ref{k-ChainRule} and \ref{t-ChainRule}). 
These chain rules will be crucial for our later arguments.

First, let us fix the following notation 
that will be used throughout this paper.

\vskip2mm

{\bf Notation:}

\vskip2mm

\begin{enumerate}

\item The base rings $k$ throughout this paper
will always be assumed to be unital 
integral domains. 

\item We fix $n\geq 1$ and let 
$z=(z_1, z_2, ... , z_n)$ be $n$ 
noncommutative variables. 
For any unital integral domain $k$,
we denote by $k\langle z\rangle $ 
and $k\langle \langle z \rangle \rangle$ the algebras of 
(noncommutative) polynomials and formal power series 
in $z_i$ $(1\leq i\leq n)$
over $k$, respectively.

\item 
For any unital integral domain $k$, 
note that the set of   
endomorphisms $\phi$ of $\kzz$ as a $k$-algebra is in 
$1$-$1$ correspondence with the set of 
$n$-vectors 
$(F_1(z), F_2(z), \cdots, F_n(z))\in \kzz^{\times n}$
via $F_i(z)=\phi(z_i)$ $(1\leq i\leq n)$. 
So, in this paper, 
by a formal endomorphism of $\kzz$ or a formal map 
in $z$, 
we simply mean a $n$-vector
$F(z)=(F_1(z), F_2(z), \cdots, F_n(z))$
with $F_i(z) \in \kzz$ $(1\leq i\leq n)$. 
When each $F_i(z)$ is a polynomial in $z$, 
we say $F(z)$ is a {\it polynomial} 
endomorphism of $\kzz$ or simply 
a {\it polynomial} map in $z$.

\item For any $m\geq 1$ and 
$U(z)=(U_1(z), \cdots, U_m(z)) \in \kzz^{\times m}$, 
we define the {\it order} $o(U(z))$ of $U(z)$ to be
\begin{align*}
o(U(z)):=\min_{1\leq i\leq k} o(U_i(z))
\end{align*}
and, when $U(z) \in \kz^{\times m}$, the {\it degree} 
$\deg (U(z))$ of $U(z)$ to be
\begin{align*}
\deg U(z):=\max_{1\leq i\leq k} \deg U_i(z).
\end{align*}
When the base ring is 
(as it frequently will be in this paper) 
the polynomial algebra $k[t]$ or 
the formal power series algebra 
$k[[t]]$ in a central parameter $t$ 
over a unital integral domain $k$,
the notation $o(U_t(z))$ and $\deg U_t(z)$ 
above always stand for the order and 
the degree of $U_t(z)$ with respect to $z$, 
respectively. In other words, 
$t$ will not be treated 
as a variable as $z_i$'s but a scalar 
parameter which commutes 
with $z_i$'s.

\item All $n$-vectors in this paper are supposed to be 
column vectors unless stated otherwise.
For any vector or matrix $U$, we denote by 
$U^\tau$ its transpose.

\end{enumerate}

Now let $k$ be a unital integral domain of any
characteristic and $\kzz$ fixed as above.   
By a {\it derivation} of $\kzz$, we mean a homomorphism 
of abelian groups 
$\delta: \kzz\to \kzz$ 
that satisfies the Leibniz rule,
i.e. for any $f, g\in \kzz$, we have
 \begin{align}\label{Leibniz}
\delta (fg)=(\delta f)g+f(\delta g).
\end{align}
A derivation $\delta$ of $\kzz$ is said to be a {\it $k$-derivation} 
if it annihilates all elements of $k\subset \kzz$. 
In other words, it is also 
a $k$-linear map from $\kzz$ to $\kzz$. 
We will denote by 
$\cDrkzz$ or $\cDrzz$, when the base ring $k$ 
is clear in the context, 
the set of all $k$-derivations of $\kzz$.
The unital subalgebra of 
$\text{End}_k(\kzz)$
generated by all
$k$-derivations 
of $\kzz$ will be denoted by $\cDzz$ or $\cDkzz$.
Elements of $\cDzz$ will be called {\it $($formal$)$ 
differential operators} in the
noncommutative variables $z_i$ 
$(1\leq i\leq n)$.

For any $1\leq i\leq n$ and $u(z)\in \kzz$, 
we denote by $\lb u(z) \fr \p{\p z_i}\rb $
the $k$-derivation which maps $z_i$ to $u(z)$ and $z_j$ to $0$ 
for any $j\neq i$. 
\footnote{The reason we put a bracket $[\cdot]$ 
in the notation for derivations of $\kzz$ is to avoid 
any possible confusion caused by a subtle point described 
in the {\bf Warning} below.}
For any $\Vec{u}=(u_1, u_2, \cdots, u_n)\in \kzz^{\times n}$, 
we set 
\begin{align}\label{Upz}
[\Vec{u}\pz]:=\sum_{i=1}^n [u_i \fr\p{\p z_i}].  
\end{align}
Furthermore, for any matrix $M_{m\times n}$ 
with row vectors 
$M_j(z)\in \kzz^{\times n}$ 
$(1\leq j\leq m)$, we set
\begin{align} \label{M-pz}
\lb M\pz \rb:=(\, [M_1 \pz], \, [M_2\pz], ... , \, [M_m\pz]\, )
 \in \cDrzz^{\times n}.
\end{align}

{\bf Warning:} {\it 
Unlike in the commutative case, 
in general, we do {\bf not} have
$\lb u(z) \fr \p{\p z_i}\rb  g(z)  = u(z)  
\fr {\p g}{\p z_i}$ for all 
$u(z), g(z)\in \kzz$.
For example, let $g = z_j z_i$ with $j\neq i$, we have
\begin{align*}
[u \fr \p{\p z_i}](z_jz_i)& = ([u \fr \p{\p z_i}]z_j )z_i+z_j 
([u \fr \p{\p z_i}]z_i)= z_j u(z),\\
u(z) \fr {\p g}{\p z_i} & =u(z)z_j,
\end{align*}
which are not equal unless $u(z)$ commutes with $z_j$. }

With the notation above, 
it is easy to see that any 
$k$-derivations $\delta$ 
of $\kzz$ can be written uniquely as 
$\sum_{i=1}^n \lb f_i(z)\fr\p{\p z_i}\rb$ 
with $f_i(z)=\delta z_i\in \kzz$ 
$(1\leq i\leq n)$.

Finally, for any automorphism $F(z)$ of $\kzz$ and
any $\delta\in \cDrzz$, we define 
$F_*(\delta)\in \cDrzz$ by setting, for any $u(z)\in \kzz$,
\begin{align}\label{star-action}
F_*(\delta)\, u(z):=\lp \delta (u(F^{-1}) \rp ( F ).
\end{align}

We call $F_*(\delta)$ the 
{\it induced action of $F(z)$ on $\delta$}.

Next, let us consider the chain rules 
for derivations of $\kzz$ and $\kttzz$.
The usual chain rule for derivations 
in the commutative case
certainly does not hold anymore 
in the noncommutative case. 
But it has the following two variations 
in certain special cases, 
see Lemma \ref{k-ChainRule} 
and  \ref{t-ChainRule} below.

First, let us consider the following chain rule 
for $k$-derivations of $\kzz$.

\begin{lemma} \label{k-ChainRule}
\, $(\text{\bf Chain Rule for $k$-Derivations})$

Let $\delta$ be a $k$-derivation of $\kzz$  
and $F(z)=(F_1(z), \cdots, F_n(z))$ 
an automorphism of $\kzz$. 
Then, for any $u(z)\in \kzz$, we have
\begin{align}\label{ChainRule1}
\delta (u(F))
= \left ( \left[ (\delta F) (F^{-1}) \fr{\p}{\p z} \right ] u \right )\circ F,
\end{align}
or equivalently, 
\begin{align}\label{ChainRule2}
(F^{-1})_*(\delta)=
 \left[ (\delta F)(F^{-1}) \fr{\p}{\p z}\right ],
\end{align}
where $\delta F:=(\delta F_1(z), \delta F_2(z), \cdots, \delta F_n(z))$.
\end{lemma}

\pf It is easy to see that Eqs.\,(\ref{ChainRule1}) 
and (\ref{ChainRule2}) are equivalent to each other 
via composing with $F$ or $F^{-1}$ from right.
So it will be enough to show Eq.\,(\ref{ChainRule2}). 

First, note that both sides of Eq.\,(\ref{ChainRule2})
are $k$-derivations of $\kzz$. Secondly, it is easy to 
check directly that,  
for any $1\leq i \leq n$, 
both derivations send $z_i$ to 
$(\delta F_i)(F^{-1})$. 
Hence they must be same as $k$-derivations
of $\kzz$ and Eq.\,(\ref{ChainRule2}) holds.
\epfv

Note that, when $z_i$'s are commutative variables, 
Eq.\,(\ref{ChainRule1}) becomes the usual chain rule.
It is worth to mention that,  
the chain rule Eq.\,(\ref{ChainRule1}) or (\ref{ChainRule2})
also has a very simple form for endomorphisms of $\kzz$ 
in terms of the Jacobian matrices.
Here,  for any sequence 
$U(z)=(U_1(z), \cdots, U_m(z))$ of $\kzz^{\times m}$, 
we define the Jacobian matrix to be
$JU(z)=\lp \lb \fr{\p }{\p z_j}\rb U_i \rp $ 
as in the commutative case and 
set $\tilde JU(z)=(JU)^\tau (z)=\lp \lb \fr{\p }{\p z_i}\rb U_j \rp $. 

\begin{corol}\label{C2.2}
Let $U(z)=(U_1, \cdots, U_m)\in \kzz^{\times m}$ and $F(z)$ 
an automorphism of $\kzz$. Then, we have
\begin{align}\label{MatrixChainRule}
\tilde J(U(F))(z)= \left( \lb \tilde JF (F^{-1}) \frac {\p}{\p z} \rb^\tau  U \right ) (F),
\end{align}
where the matrix $\left( \lb \tilde JF (F^{-1}) \frac {\p}{\p z} \rb^\tau  U \right )$
in the equation above is the formal ``product" of the column 
vector $\lb \tilde JF (F^{-1}) \frac {\p}{\p z} \rb^\tau \in  \cDrzz^{\times n}$ 
with the row vector $U(z)= (U_1, \cdots, U_m)$.

In particular, when $m=n$ and $U(z)=G(z):=F^{-1}(z)$, 
 we have 
\begin{align}\label{JacobianChainRule}
\left[ \tilde JF(G) \frac {\p}{\p z}\right] G =I_{n\times n}
=\left[ \tilde JG (F) \frac {\p}{\p z} \right] F(z).
\end{align}
\end{corol}

The proof of Eq.\,(\ref{MatrixChainRule}) 
is straightforward, just to apply 
Eq.\,(\ref{ChainRule1}) or (\ref{ChainRule2})
to the entries of the matrix 
$\tilde J(U(F))(z)$; while 
Eq.\,(\ref{JacobianChainRule}) 
is an immediate consequence of 
Eq.\,(\ref{MatrixChainRule}) 
and the fact $G(F(z))=z=F(G(z))$.

Note that, when $z$ are commutative variables, 
Eq.\,(\ref{JacobianChainRule}) is same as
$JF(G)JG=I_{n\times n}=JG(F)JF$.
But, unlike in the commutative case, 
$JF(G)$ in general
is not the multiplication inverse matrix of $JF$. 
This can be seen from the following example.

\begin{exam}
Let $F(x, y)=(F_1, F_2)$ be the automorphism of $k\langle\langle x, y\rangle\rangle$ 
with
\begin{align*}
F_1(x, y)&= e^x-1,\\
F_2(x, y)&= ye^{-x}.
\end{align*}
Its inverse map $G(x, y)=(G_1, G_2)$ is given by 
\begin{align*}
G_1(x, y)&= \ln(1+x),\\
G_2(x, y)&= y(1+x).
\end{align*}
Now consider the Jacobian matrices 
\begin{align*}
JF(x, y)=
\begin{pmatrix} 
e^x & 0\\
-ye^{-x} & e^{-x}
\end{pmatrix},
\quad\quad
JG(x, y)=
\begin{pmatrix} 
\fr 1{1+x} & 0\\
y & 1+x
\end{pmatrix}
\end{align*}
But, on the other hand, 
\begin{align*}
JG(F_1, F_2)=
\begin{pmatrix} 
e^{-x} & 0\\
ye^{-x} & e^{x}
\end{pmatrix},
\quad \quad
(JF)^{-1}(x, y)=
\begin{pmatrix} 
e^{-x} & 0\\
e^xye^{-2x} & e^{x}
\end{pmatrix}
\end{align*}

Hence $JG(F)\neq (JF)^{-1}$ unless $x$ and $y$ 
commute with each other. 
\end{exam}

The second chain rule we will need later is the following.
Let $t$ be a formal parameter which commutes with $z$ and 
$k[[t]]$ the formal power series in $t$ over $k$. Note that
the derivation $\fr\p{\p t}$ of $k[[t]]$ can be
extended naturally 
to a derivation 
of $\kttzz$, which we still denote by $\fr\p{\p t}$,  
by setting $\fr{\p z_i}{\p t}=0$ 
for any $1\leq i\leq n$. 

\begin{lemma}\label{t-ChainRule}
Let $F_t=(F_{t, 1}, F_{t, 2},\cdots, F_{t, n} )$ be
an automorphism of $\kttzz$ $(\text{as an algebra over 
$k[[t]]$})$ with inverse map $F^{-1}_t(z)$.
Then, for any $u_t(z)\in \kttzz$, we have

\begin{align}\label{t-ChainRule-e1}
\frac {\p u_t(F_t)}{\p t}
&=\frac {\p u_t}{\p t}(F_t)+
\left( \left[ \fr{\p  F_{t}}{\p t} 
(F_t^{-1}) \fr{\p}{\p z}   
\right ] u_t \right )(F_t).
\end{align}
\end{lemma}

\pf The proof is similar as the one for Lemma \ref{k-ChainRule}, 
which goes as follows. 

First, composing $F_t^{-1}$ to Eq.\,(\ref{t-ChainRule-e1}) from right, 
we get
 \begin{align}\label{t-ChainRule-F}
\frac {\p u_t(F_t)}{\p t}\circ F_t^{-1}
&=\frac {\p u_t}{\p t}(z)+
 \left[ \fr{\p  F_{t}}{\p t} 
(F_t^{-1}) \fr{\p}{\p z}   
\right ] u_t, 
\end{align}
which is equivalent to Eq.\,(\ref{t-ChainRule-e1}).
 
Secondly, we define the maps 
$\delta_1, \delta_2:\kttzz \to \kttzz$ by setting 
\begin{align}
\delta_1(u_t)&=\frac {\p u_t(F_t)}{\p t}\circ F_t^{-1}, \\
\delta_2(u_t)&=\frac {\p u_t}{\p t} +
 \left[ \fr{\p  F_{t}}{\p t} 
(F_t^{-1}) \fr{\p}{\p z}   
\right ] u_t \end{align}
for any $u_t(z) \in \kttzz$. 

Hence, it will be enough to show $\delta_1=\delta_2$. 
But, again, it is easy to see  
that $\delta_i$ $(i=1, 2)$ both
are derivations of $\kttzz$. 
(Actually, $\delta_1=(F_t^{-1})_*(\frac\p{\p t})$). 
Therefore, 
it will be enough to show 
they have same values 
when $u_t(z)=t$ and $u_t(z)=z_i$ 
for any $1\leq i\leq n$. 
But, for these cases, we have
\begin{align*}
\delta_1(t)&=1=\delta_2(t),\\
\delta_1(z_i)&= \fr{\p F_{t,i}}{\p t}(F_t^{-1})=\delta_2(z_i)
\end{align*}
for any $1\leq i\leq n$.
\epfv

\renewcommand{\theequation}{\thesection.\arabic{equation}}
\renewcommand{\therema}{\thesection.\arabic{rema}}
\setcounter{equation}{0}
\setcounter{rema}{0}

\section{\bf General Deformations} \label{S3}

Let $k$ be 
a unital integral domain  
of any characteristic and 
$z=(z_1, z_2, ... , z_n)$ 
and $t$ as in the previous section, 
i.e. $z_i$ $(1\leq i\leq n)$ are
$n$ free noncommutative 
variables and 
$t$ is a formal 
parameter which commutes 
with $z_i$'s. 
In this section,
we study the general 
deformation of 
automorphisms of $\kzz$ 
parameterized by $t$.
More precisely, 
we study automorphisms
$F_t(z)$ of $\kttzz$ over $k[[t]]$
of the form $F_t(z)=z-H_t(z)$ with 
$H_t(z)\in \kttzz^{\times n}$ and 
$o(H_t(z))\geq 2$.
Note that, when $F(z):=F_{t=1}(z)$
makes sense 
(for example, when $H_t(z)\in \ktzz^{\times n}$, or 
when $k=\bC$ and all coefficients of $H_t(z)$ 
are holomorphic functions of $t$ which are 
convergent over an open subset of $\bC$ 
containing the closed unit disk),
$F_t(z)$ can be viewed 
as a deformation of the automorphism 
$F(z):=F_{t=1}(z)$ of $\kzz$. 
We will denote by $G_t(z)$ 
and $G(z)$ the formal inverse maps 
of $F_t(z)$ and 
$F(z)=F_{t=1}(z)$ 
(again, when it exists), 
respectively.  
We will always write $G_t(z)$ 
as $G_t(z)=z+M_t(z)$ for some 
$M_t(z)\in \ktz^{\times n}$ 
with $\mbox{o} (M_t(z))\geq 2$. 
Note that, when $F(z)=F_{t=1}(z)$ 
and $G_{t=1}(z)$ both make sense, 
by the uniqueness 
of inverse maps, 
we have $G_{t=1}(z)=G(z)$. 
In this section,
we first derive  
the PDE's satisfied by $M_t(z)$, 
$u(F_t)$ and $u(G_t)$ with 
$u(z)\in \kzz$ 
(see Eqs.\,(\ref{L3.1.1-e4}), (\ref{P3.1.3-e1}) 
and (\ref{P3.1.3-e2})).
We then in Theorem \ref{T3.1.4} show that,  
when $char.\, k=0$, 
the power series $u(F_t)$ and $u(G_t)$ 
$(u(z)\in \kzz)$ are actually 
characterized 
by the PDE's (\ref{P3.1.3-e1}) and 
(\ref{P3.1.3-e2}), respectively. 
When $char.\, k=p>0$, 
$u(F_t)$ and $u(G_t)$ still satisfy
the PDE's (\ref{P3.1.3-e1}) and 
(\ref{P3.1.3-e2}), respectively
but they are only 
uniquely determined 
by these PDE's together
with their coefficients
of $t^{mp}$ $(m\geq 0)$
(see Remark \ref{R3.1.5}).

\vskip2mm

Let us start with the following simple lemma.

\begin{lemma}\label{L3.1.1}
Let $F_t(z)$, $H_t(z)$, $G_t(z)$, $M_t(z)$ 
as fixed above. Then we have
\allowdisplaybreaks{
\begin{align}
M_t&=H_t(G_t),\label{L3.1.1-e1} \\
H_t&=M_t(F_t),\label{L3.1.1-e2} \\
\frac{\p H_t}{\p t}(z) &=
\left[\frac{\p M_t}{\p t}(F_t) \pz \right ] F_t(z) , \label{L3.1.1-e3}\\
\frac{\p M_t}{\p t}(z) &=
\left[\frac{\p H_t}{\p t}(G_t) \pz \right ] G_t(z).\label{L3.1.1-e4}
\end{align}}
\end{lemma}
\pf
Since $F_t(G_t(z))=z$, we have 
\begin{align}
z+M_t(z)-H_t(G_t(z))=z.
\end{align}
Hence Eq.\,$(\ref{L3.1.1-e1})$ holds. Similarly, Eq.\,$(\ref{L3.1.1-e2})$ 
follows from $G_t(F_t(z))=z$.

To show Eq.\,$(\ref{L3.1.1-e3})$, applying $\frac{\p}{\p t}$ to Eq.\,$(\ref{L3.1.1-e1})$ 
and using the chain rule Eq.\,$(\ref{t-ChainRule-e1})$, 
we have 
\begin{align*}
\frac{\p M_t}{\p t}&= \frac{\p H_t(G_t)}{\p t}\\
&=\frac{\p H_t}{\p t}(G_t)+
\left( \left[ \fr{\p  G_{t}}{\p t}(F_t) \fr{\p}{\p z}   
\right ] H_t \right )(G_t), \\
&= \frac{\p H_t}{\p t}(G_t)+
\left( \left[ \fr{\p  M_{t}}{\p t}(F_t) \fr{\p}{\p z}   
\right ] H_t \right )(G_t).
\end{align*}
Therefore, we have 
\begin{align*}
\frac{\p H_t}{\p t}(G_t)&=
\frac{\p M_t}{\p t}-
\left( \left[ \fr{\p  M_{t}}{\p t}(F_t) \fr{\p}{\p z}   
\right ] H_t \right )(G_t)\\
&=\left( \left[ \fr{\p  M_{t}}{\p t}(F_t) \fr{\p}{\p z}   
\right ] (z-H_t) \right )(G_t) \\
&=\left( \left[ \fr{\p  M_{t}}{\p t}(F_t) \fr{\p}{\p z}   
\right ] F_t \right )(G_t).
\end{align*}

Composing with $F_t$ from right to the equation above, we get Eq.\,$(\ref{L3.1.1-e3})$.
Eq.\,$(\ref{L3.1.1-e4})$ can be proved similarly by applying $\frac{\p}{\p t}$ to 
Eq.\,$(\ref{L3.1.1-e2})$. 
\epfv 

Now, we set 
\begin{align}
h(t)&:=\left[\frac{\p M_t}{\p t}(F_t) \pz \right ],\label{Def-h(t)}\\
m(t)&:=\left[\frac{\p H_t}{\p t}(G_t) \pz \right ].\label{Def-m(t)}
\end{align}

\begin{lemma}\label{L3.1.2}
\begin{align}
(G_t)_* (h(t))&=m(t),\label{L3.1.2-e1}\\
(F_t)_* (m(t))&=h(t).\label{L3.1.2-e2}
\end{align}
\end{lemma}
\pf Note that Eq.\,$(\ref{L3.1.2-e2})$ follows immediately 
when we apply $(F_t)_*$ to Eq.\,$(\ref{L3.1.2-e1})$.
So we only need show Eq.\,$(\ref{L3.1.2-e1})$.

First, applying the chain rule 
Eq.\,$(\ref{ChainRule2})$ 
with $\delta=h(t)$ and Eq.\,$(\ref{Def-h(t)})$, 
we have 
\allowdisplaybreaks{
\begin{align*}
(G_t)_* (h(t))&= \left[ \lp h(t)F_t \rp (G_t) \pz \right ] \\
&=\left[ \lp 
\left[\frac{\p M_t}{\p t}(F_t) \pz \right ] F_t \rp
(G_t) \pz \right ] \\
\intertext{Applying Eqs.\,$(\ref{L3.1.1-e3})$ and $(\ref{Def-m(t)})$:}
&= \left[\frac{\p H_t}{\p t}(G_t) \pz \right ] \\
&=m(t).
\end{align*} }
\epfv

\begin{propo}\label{P3.1.3}
For any $u(z)\in \kzz$, we have
\begin{align}
\frac {\p \, u(F_t) }{\p t}  
&=-(m(t)u)(F_t)=-h(t)\,\, u(F_t),\label{P3.1.3-e1}\\
\frac {\p \, u(G_t )}{\p t} 
&=(h(t)u)(G_t)=m(t)\,\, u(G_t).\label{P3.1.3-e2}
\end{align}
\end{propo}

\pf Here we only give a proof for Eq.\,$(\ref{P3.1.3-e1})$. 
Eq.\,$(\ref{P3.1.3-e2})$ can be proved 
by a similar argument. 

By the chain rule Eq.\,$(\ref{t-ChainRule-e1})$, we have

\begin{align}\label{P3.1.3-pe1}
\frac {\p \, u(F_t)}{\p t}
&=\frac {\p u}{\p t}(F_t)+
\left( \left[ \fr{\p  F_{t}}{\p t} 
(G_t) \fr{\p}{\p z}   
\right ] u \right )(F_t) \\
&=-\left( \left[ \fr{\p  H_{t}}{\p t} 
(G_t) \fr{\p}{\p z}   
\right ] u \right )(F_t)\nno \\
&=-(m(t)u)( F_t). \nno
\end{align}
Hence, we get the first part of 
Eq.\,$(\ref{P3.1.3-e1})$. 
To show the second part, 
first, by Eqs.\,$(\ref{L3.1.2-e1})$ and $(\ref{star-action})$, 
we have
\begin{align*} 
 m(t)u(z)
&=\lp (G_t)_* h(t)\rp u(z) \\
&=\left ( h(t)\,\, u(F_t)\right)(G_t).\nno
\end{align*}
Composing with $F_t$ from right 
to the equation above, we get
\begin{align} \label{P3.1.3-pe2}
 (m(t)u)(F_t)= h(t)\,\, u(F_t).
\end{align}
Combining Eqs.\,(\ref{P3.1.3-pe1}) and (\ref{P3.1.3-pe2}), 
we have
\begin{align*}
\frac {\p u(F_t)}{\p t}
&=-(m(t)u)( F_t)\\
&=-h(t)\,\, u(F_t),
\end{align*}
which is the second part of Eq.\,$(\ref{P3.1.3-e1})$. 
\epfv 

Actually, when char.\,$k=0$, 
elements of $\kttzz$ of the forms $u(F_t)$ 
and $u(G_t)$ for some $u(z)\in \kzz$
are characterized by 
Eqs.\,$(\ref{P3.1.3-e1})$ and 
$(\ref{P3.1.3-e2})$, respectively. 
This can be seen from the following theorem.

\begin{theo}\label{T3.1.4}
Assume that the base ring $k$ has char.\,$k=0$, then

$(a)$ For any $U_t(z)\in \kttzz$, $U_t(z)=u(F_t(z))$ for some 
$u(z)\in \kzz$ iff $U_t(z)$ satisfies the PDE
\begin{align}
\frac {\p U_t(z)}{\p t}  
&=  -h(t) U_t(z). \label{T3.1.4-e1}
\end{align}

$(b)$ For any $V_t(z)\in \kttzz$, $V_t(z)=u(G_t(z))$ for some 
$u(z)\in \kzz$ iff $V_t(z)$ satisfies the PDE
\begin{align}
\frac {\p V_t(z)}{\p t}  
&=  m(t) V_t(z). \label{T3.1.4-e2}
\end{align}
\end{theo}

\pf $(a)$ The $(\Rightarrow)$ part is just Proposition \ref{P3.1.3}. 
Conversely,
suppose $U_t(z)\in \kttzz$ satisfies Eq.\,$(\ref{T3.1.4-e1})$. 
Set $\wtilde U_t(z)=U_t(G_t(z))$.
By the chain rule Eq.\,$(\ref{t-ChainRule-e1})$, we have
\allowdisplaybreaks{
\begin{align*}
\frac {\p \wtilde U_t(z)}{\p t}&= 
\frac {\p U_t}{\p t}(G_t)+ \lp \lb \frac{\p G_t}{\p t}(F_t) \pz \rb U_t \rp (G_t)\\
&= \left( \frac {\p U_t}{\p t} + \lb \frac{\p G_t}{\p t}(F_t) \pz \rb U_t \right) (G_t)\\
&= \left( \frac {\p U_t}{\p t} + h(t)U_t  \right) (G_t)\\
&=0. 
\end{align*}}
Therefore, if we set $u(z):=\wtilde U_t(z)=U_t(G_t(z))$, 
then $u(z)\in \kzz$ and $U_t(z)=u(F_t)$. 
Hence we have proved $(a)$. 

$(b)$ can be proved similarly.
\epfv

\begin{rmk}\label{R3.1.5} 
From the proof of Theorem \ref{T3.1.4} above, 
one can see that,
when the base ring $k$ has $char.\,k=p>0$, 
the $(\Rightarrow)$ part of 
the theorem still holds; while  
the $(\Leftarrow)$ part 
is not true in general.
But, if the coefficients 
of $t^{mp}$ $(m\geq 0)$ of 
$U_t(z)$ and $V_t(z)$ 
are given or fixed, 
$U_t(z)$ and $V_t(z)$ are still 
uniquely determined by 
Eqs.\,$(\ref{T3.1.4-e1})$ and 
$(\ref{T3.1.4-e2})$, 
respectively. 
This can be easily seen by 
viewing $U_t(z)$ and $V_t(z)$ 
as formal power series in $t$ over the ring $\kzz$ 
and solving  Eqs.\,$(\ref{T3.1.4-e1})$ and 
$(\ref{T3.1.4-e2})$ recursively. For a more 
detailed discussion on 
a similarly situation, see Section \ref{S5}.
\end{rmk}

\renewcommand{\theequation}{\thesection.\arabic{equation}}
\renewcommand{\therema}{\thesection.\arabic{rema}}
\setcounter{equation}{0}
\setcounter{rema}{0}

\section{\bf A Special Deformation} \label{S4}

In this section, we will focus on 
a special family of deformations
of automorphisms of $\kzz$. We start with a fixed
automorphism $F(z)$ of $\kzz$ 
and always assume that $F(z)$ 
has the form $F(z)=z-H(z)$ with 
$o(H(z))\geq 2$. We set 
$F_t(z)=z-tH(z)$ and write 
its inverse map as $G_t(z)=z+tN_t(z)$ 
with $N_t(z)\in \kttzz^{\times n}$ and $o(N_t(z))\geq 2$. 
In terms of the notation in Section \ref{S2}, 
we have
\begin{align}
H_t(z)&=tH(z), \label{tH(z)}\\
M_t(z)&=tN_t(z).\label{tN_t(z)}
\end{align}

We first apply the results obtained 
in the previous section for 
the general deformations to 
the special deformation above.
In particular, we show in Theorem \ref{T4.3} 
that $N_t(z)$ is a power series 
solution of a Cauchy problem 
of the PDE involved 
(see Eqs.\,(\ref{NC-PDE}) and (\ref{NC-PDE-B})). 
One interesting aspect of 
this fact is that, 
when passing to the 
commutative case, the PDE 
(\ref{NC-PDE}) is almost the Burgers' equation 
in Diffusion theory, see Remark \ref{rmk-Burgers}. 
When $F_t(z)=z-tH(z)$ is a symmetric map, 
i.e. $H(z)$ is the gradient vector $\nabla P(z)$ 
for some $P(z)\in k[[z]]$,
it can be further linked to the Heat equation. 
For more discussion in this direction, 
see \cite{BurgersEq} and \cite{HNP}.
The PDE (\ref{NC-PDE}) 
in Theorem \ref{T4.3} is also
the starting point for the  
inversion formulas 
that will be derived 
in next two sections.
Besides the property of $N_t(z)$ given 
in Theorem \ref{T4.3}, 
other characterizing properties of $N_t(z)$ 
are also derived 
(see Lemma \ref{L4.6} and Proposition \ref{P4.7}).

\vskip2mm

First, let us work out the special 
forms for the differential operators 
$h(t)$ and $m(t)$ defined in 
Eqs.\,$(\ref{Def-h(t)})$ and 
$(\ref{Def-m(t)})$, respectively, 
for the special deformation 
$F_t(z)=z-tH(z)$ with 
$H(z)\in \kzz^{\times n}$ 
and $o(H(z))\geq 2$.

\begin{lemma}\label{L4.1}
With the notation above, we have
\begin{align}
m(t)&=\lb N_t(z)\pz \rb, \label{L4.1-e1} \\
h(t)&=\sum_{m\geq 1}  t^{m-1} \lb C_m(z)\pz \rb, \label{L4.1-e2}
\end{align}
where $C_m(z)\in \kzz^{\times n}$ $(m\geq 1)$ are defined recursively by
\begin{align}
C_1(z)& =H(z),\label{L4.1-e3}\\
C_m(z)&=\lb C_{m-1}(z)\pz \rb H,\label{L4.1-e4}
\end{align}
for any $m\geq 2$.
\end{lemma}

\pf First, by Lemma \ref{L3.1.1} and Eqs.\,(\ref{tH(z)}) and (\ref{tN_t(z)}), 
it is easy to see that, we have
\BQ
N_t(F_t(z))&=&H(z), \label{L4.1-pe1}\\
H(G_t)&=& N_t(z). \label{L4.1-pe2}
\EQ

By Eqs.\,(\ref{Def-m(t)}), (\ref{tH(z)}) and also the equations above,
we have
\begin{align*}
m(t)&=\left[\frac{\p H_t}{\p t}(G_t) \pz \right ]\\
&=\left[ H(G_t) \pz \right ]\\
&=\lb N_t(z)\pz \rb. 
\end{align*}
Hence, we get Eq.\,(\ref{L4.1-e1}).

To show Eq.\,(\ref{L4.1-e2}), we first write $h(t)$ as in 
Eq.\,$(\ref{L4.1-e2})$ for some $C_m(z)\in \kzz^{\times n}$ 
$(m\geq 1)$,  and then show that $C_m(z)$'s also satisfy
Eqs.\,(\ref{L4.1-e3}) and (\ref{L4.1-e4}). Consequently,
$C_m(z)$ $(m\geq 1)$
will be uniquely determined 
by Eqs.\,(\ref{L4.1-e3}) 
and (\ref{L4.1-e4}).

First, by Eqs.\,(\ref{L3.1.1-e3}) and (\ref{Def-h(t)}), we have
\allowdisplaybreaks{
 \begin{align*}
H(z)&=\frac{\p H_t}{\p t}(z)\\
&=h(t)F_t(z) \\
&= \sum_{m\geq 1}  t^{m-1} \lb C_m(z) \pz \rb (z-tH(z))\\
&=\sum_{m\geq 1}  t^{m-1} C_m(z)-t\sum_{m\geq 1}  t^{m-1} \lb C_m(z) \pz \rb H(z)\\
 &=C_1(z)+ \sum_{m\geq 2}  t^{m-1} \lp C_m(z)-\lb C_{m-1}(z) \pz \rb H(z) \rp. 
\end{align*} }
Then, by comparing the coefficients of $t^{m-1}$ 
$(m\geq 1)$ in the equation above, we see that $C_m(z)$ 
$(m\geq 1)$ indeed satisfy Eqs.\,(\ref{L4.1-e3}) 
and (\ref{L4.1-e4}).
\epfv

By using the mathematical 
induction on $m\geq 1$, 
it is easy to check that, 
when $z_i$'s are commutative variables, 
$C_m(z)$ further has 
the following simple form.
 
\begin{corol}\label{C4.2}
For commutative variables $z_i$ $(1\leq i\leq n)$, we have
\begin{align}
C_m(z)=(JH)^{m-1}H,\label{C4.2-e1}
\end{align}
for any $m\geq 1$.
\end{corol}

By Eqs.\,(\ref{L4.1-e1}), (\ref{L4.1-pe2}) and 
Theorem \ref{T3.1.4}, 
$(b)$ with $u(z)=H_i(z)$ $(1\leq i\leq n)$ 
for the special deformation $F_t$, 
it is easy to see that we have 
the following theorem on $N_t(z)$, 
which later will imply 
an effective recurrent 
inversion formula for $G_t(z)$ 
(see Theorem \ref{T5.5}).

\begin{theo}\label{T4.3}
Let $k$ be a unital integral domain of any characteristic  
and $H(z)\in \kzz ^{\times n}$, 
$N_t(z)\in \kttzz^{\times n}$ 
as above, then, 
$N_t(z)$ is a power series solution in $\kttzz^{\times n}$ of the following 
Cauchy problem of PDE's in noncommutative variables.
\BQ
 &{}& \frac {\p N_t}{\p t}= \left[ N_t \pz \right] N_t \label{NC-PDE}\\
 &{}& N_{t=0}(z)=H(z). \label{NC-PDE-B}
\EQ
\end{theo}

\begin{rmk}\label{rmk-Burgers}
Note that, in the commutative case, 
the PDE $(\ref{NC-PDE})$ 
becomes
\begin{align}
 \frac {\p N_t}{\p t}= J N_t\cdot N_t \label{COM-PDE}.
\end{align}
which was first proved in \cite{Recurrent} 
$(\text{unpublished})$ 
and later in \cite{BurgersEq}. Interestingly, 
the PDE above is almost the  classical Burgers' 
equation in Diffusion theory, 
which has the form
\BQ
 &{}& \frac {\p N_t}{\p t}= \lp J N_t\rp^\tau \cdot N_t \label{Burgers}.
\EQ
In particular, when $N_t$ is the gradient vector of $Q_t$ for some
$Q_t\in k[[t]][[z]]$, Eqs.\,$(\ref{COM-PDE})$  
and $(\ref{Burgers})$ coincide. 
Furthermore, 
in this case,  Eq.\,$(\ref{COM-PDE})$
is also closely related with the Heat equation.
For more detailed discussions on the connections 
among these three PDE's in the commutative case,  
see \cite{BurgersEq} and \cite{HNP}. 
\end{rmk}

Next, we derive more properties of 
$N_t(z)$. The first interesting property 
of $N_t(z)$ is the following proposition.
It essentially says that $\{N_t(z)| t\in k \}$ gives 
a family of automorphisms of $\kttzz$ 
which are ``closed'' under 
the inverse operation.

\begin{propo}\label{P4.4}
For any $s\in k$, 
 the formal inverse of $U_{s, t}(z):=z-sN_t(z)$
is given by
$V_{s, t}(z):=z+ s N_{t+s}(z)$.
Actually, $U_{s, t}(z)= F_{t+s}\circ G_t(z)$ 
and $V_{s, t}(z)= F_{t}\circ G_{s+t}(z)$.
\end{propo}
\pf 
\BQn
F_{t+s}\circ G_t (z)&=& G_t(z)-(t+s)H(G_t(z))\\
&=& z+tN_t(z)-(t+s)N_t(z)\\
&=& z-sN_t(z)\\
&=&U_{s, t}(z).
\EQn
 Similarly, we can prove $V_{s, t}(z)= F_{t}\circ G_{s+t}(z)$. Hence we have
$U_{s, t}^{-1}(z)=V_{s, t}(z)$.
\epfv

In the rest of this section, we 
will assume the base ring $k$ 
has char.\,$k=0$.
Below we show that $N_t(z)$ in this case
is actually characterized by 
the Cauchy problem Eqs.\,$(\ref{NC-PDE})$ 
and $(\ref{NC-PDE-B})$ in Theorem \ref{T4.3}.

\begin{propo}\label{P4.5}
For any $H(z)\in \kzz ^{\times n}$ and $N_t(z)\in \kttzz^{\times n}$ with
$\mbox{o} (H(z))\geq 2$ and $\mbox{o} (N_t(z))\geq 2$, respectively. 
The following statements are equivalent.
\begin{enumerate}
\item The formal map $G_t(z)=z+tN_t(z)$ is the inverse of
$F_t(z)=z-tH(z)$.

\item $N_t(z)\in \kttzz$ is the unique power series solution of 
the Cauchy problem Eqs.\,$(\ref{NC-PDE})$ and $(\ref{NC-PDE-B})$.
\end{enumerate}
\end{propo}

\pf First, $(1)\Rightarrow (2)$ is exactly Theorem \ref{T4.3}. 
To show $(2)\Rightarrow (1)$, we assume that the formal 
inverse of $F_t(z)=z-tH(z)$ is given by 
$G_t(z)=z+t\widetilde N_t(z)$. By Theorem \ref{T4.3}, 
we know that 
$\widetilde N_t(z)$ also satisfies Eqs.\,(\ref{NC-PDE}) 
and (\ref{NC-PDE-B}). But, by Corollary \ref{Uniqueness}, 
$(a)$ in next section, 
the power series solution in $\kttzz$ of
Eqs.\,(\ref{NC-PDE}) and (\ref{NC-PDE-B}) 
is actually unique. 
Hence we have $\wtilde N_t(z)=N_t(z)$ and  
$(2)\Rightarrow (1)$ follows.
\epfv

Another characterizing property of $N_t(z)$ 
(see Proposition \ref{P4.7} below) can be derived 
as follows. First, we need the following lemma.

\begin{lemma}\label{L4.6}
For any $u(z)\in \kzz $, the unique power series 
solution $U_t(z)$ in $z$ and $t$ of the following 
Cauchy problem

\begin{align}\label{GPDE}
\begin{cases}
&\frac {\p U_t}{\p t} = \left[ N_t \pz \right] U_t, \\
&U_{t=0}(z) = u(z).
\end{cases}
\end{align} 
is given by $U_t(z)=u(z+tN_t(z))$.
\end{lemma}

\pf By a similar argument as in
the proof of Lemma \ref{L5.1} in next section, 
it is easy to check that the power series solution 
in $z$ and $t$ of the Cauchy problem 
Eq.\,(\ref{GPDE}) is unique. 
So it will be enough to 
show that $U_t(z)=u(z+tN_t(z))$ 
satisfies Eq.\,(\ref{GPDE}). 
First, the boundary 
condition in 
Eq.\,(\ref{GPDE}) 
is obviously 
satisfied by $U_t(z)$.
Secondly, by Theorem \ref{T3.1.4}, $(b)$ 
and Eq.\,$(\ref{L4.1-e1})$, 
$U_t(z)$ also satisfies 
the PDE in Eq.\,(\ref{GPDE}). 
\epfv

\begin{propo}\label{P4.7}
For any $N_t(z)\in \kttzz^{\times n}$ with $o(N_t(z))\geq 2$, the following 
two statements are equivalent.

$(a)$ $z+tN_t(z)$ is the formal inverse map 
of $F_t(z)\!=\!z-tH(z)$ for some $H(z)\in \kzz^{\times n}$.

$(b)$ Lemma \ref{L4.6} holds for $N_t(z)$.
\end{propo}

\pf First, $(a)\Rightarrow (b)$ follows from 
Lemma \ref{L4.6}. 
To show $(b)\Rightarrow (a)$,   
let $U_{t, i}(z)$ $(1\leq i\leq n)$ be the unique 
power series
solution of the Cauchy problem (\ref{GPDE})
with $u(z)=z_i$. Set 
$\wtilde U_t(z)=(U_{t, 1}(z), \cdots, U_{t, n}(z))$.
Note that Eq.\,(\ref{GPDE}) for $U_{t, i}(z)$ $(1\leq i\leq n)$ 
can be written as
\begin{align}
\frac {\p \wtilde U_t}{\p t} & = \left[ N_t \pz \right] \wtilde U_t. \label{GPDE-2}
\end{align}

Since, by our condition on $N_t(z)$, Lemma \ref{L4.6} holds for $N_t(z)$, 
so we have
\begin{align}\label{GPDE-3}
\wtilde U_{t}(z)&=z + tN_{t}(z).
\end{align}

Applying $\fr {\p}{\p t}$ to the equation above, we get
\BQ
\fr {\p \wtilde U_t}{\p t}=N_t+t \fr {\p N_t}{\p t}.
\EQ
Combining the equation above with Eqs.\,(\ref{GPDE-2}) and (\ref{GPDE-3}), 
we have 
\begin{align*}
N_t+t \fr {\p N_t}{\p t}= \left[ N_t \pz \right] (z + tN_{t}) 
=N_t+t\left[ N_t \pz \right] N_t.
\end{align*}
Therefore, we have
\begin{align}
\fr {\p N_t}{\p t}=\left[ N_t \pz \right] N_t.
\end{align}
Set $H(z)\!:=\! N_{t=0}(z)$. Therefore, $N_t(z)$ is a 
formal power series solution of the Cauchy problem 
Eqs.\,$(\ref{NC-PDE})$ and $(\ref{NC-PDE-B})$. 
Then, by Proposition 
\ref{P4.5}, we see that
$(a)$ holds.
\epfv

\renewcommand{\theequation}{\thesection.\arabic{equation}}
\renewcommand{\therema}{\thesection.\arabic{rema}}
\setcounter{equation}{0}
\setcounter{rema}{0}

\section{\bf A Recurrent Inversion Formula 
for automorphisms in Noncommutative Variables } \label{S5}

In this section, we apply 
some results obtained in 
Section \ref{S4} 
to derive a recurrent 
inversion formula for formal maps 
in noncommutative variables 
(see Theorem \ref{T5.5}). 
This will generalize 
the recurrent inversion 
formula in \cite{BurgersEq} 
for the commutative case with 
$char.\,k=0$ to the noncommutative case 
over a base ring $k$ of any characteristic.

\begin{lemma} \label{L5.1}
Let $W_t(z)\in \kttzz$ be a solution  
of Eqs.\,$(\ref{NC-PDE})$ and $(\ref{NC-PDE-B})$. 
We write $W_t(z)$ as
\begin{align} \label{L5.1-e1}
W_t(z)= \sum_{m=1}^\infty  W_{[m]}(z) t^{m-1}.
\end{align}
with $W_{[m]}(z)\in k\langle \langle z\rangle\rangle$ $(m\geq 1)$.
Then, the sequence $\{W_{[m]}(z) | m\geq 1\}$
satisfies the following recurrent relations: 
\BQ
W_{[1]}(z) &=& H(z), \label{M1}\\
(m-1) W_{[m]}(z) &=&  \sum_{\substack{k+l=m\\ k, l\geq 1}} 
 \left [ W_{[k]}\pz \right ] W_{[l]}  \label{Mm}
\EQ
for any $m\geq 2$.
\end{lemma}

\pf First, Eq.\,(\ref{M1}) follows directly 
from Eq.\,$(\ref{NC-PDE-B})$. 
Secondly, by Eq.\,(\ref{NC-PDE}), we have
\begin{align*}
\sum_{m=1}^\infty (m-1) W_{[m]}(z)t^{m-2} =
\left (\sum_{k=1}^\infty  t^{k-1} \left [ W_{[k]} \pz \right ] \right )
\left (\sum_{l=1}^\infty  W_{[l]}(z) t^{l-1} \right ).
\end{align*}
For any $m\geq 2$, by
comparing the coefficients of $t^{m-2}$ of the both sides 
of the equation above, we get Eq.\,(\ref{Mm}). 
\epfv

Some direct consequences 
of the lemma above are given by
the following three corollaries.

\begin{corol}\label{Uniqueness}
$(a)$ When char.\,$k=0$, the power series solutions in 
$\kttzz$ of the Cauchy problem Eqs.\,$(\ref{NC-PDE})$ and 
$(\ref{NC-PDE-B})$ 
is unique. 

$(b)$ When char.\,$k=p>0$, there are infinitely many 
solutions $W_t(z)$ in $\kttzz$ of the Cauchy problem 
Eqs.\,$(\ref{NC-PDE})$ and $(\ref{NC-PDE-B})$. 
Actually, for any fixed 
$W_{[mp+1]}\in \kzz$ $(m\geq 1)$, there exists one and only one 
solution of Eqs.\,$(\ref{NC-PDE})$ and $(\ref{NC-PDE-B})$.
\end{corol}

Let $H(z)$ and $N_t(z)$ be fixed as 
in Section \ref{S4}. 
We define the sequence 
$\{ N_{[m]}(z)\in \kzz |m\geq 1\}$ by writing 
\begin{align}\label{Def-Nm}
N_t(z)=\sum_{m\geq 1} t^{m-1} N_{[m]}(z).
\end{align}

\begin{corol} \label{C5.3}
Suppose that the base ring $k$ has char.\,$k=p>0$. Then, 
for any $m\geq 1$ and $m\equiv 1\mod p$, 
we have
\begin{align} \label{C5.3-e1}
 \sum_{\substack{k+l=m\\ k, l\geq 1}}  
\left [ N_{[k]}\pz \right ] N_{[l]}(z)=0.
\end{align}
\end{corol}
\pf 
By Theorem \ref{T4.3} 
and Lemma \ref{L5.1}, we know
the sequence 
$\{ N_{[m]}(z)\in \kzz |m\geq 1\}$
satisfies the recurrent relations 
Eqs.\,(\ref{M1}) and (\ref{Mm}).
Hence the corollary follows immediately from 
Eq.\,(\ref{Mm}).
\epfv

\begin{corol}\label{C5.4}
For any unital integral domain $k$ of any 
characteristic, we have

$(a)$ $o(N_{[m]}(z)) \geq m+1$ for any $m\geq 1$.

$(b)$ Suppose $H(z)\in \kz^{\times n}$, then , for any $m\geq 1$, 
 $N_{[m]}(z) \in \kz ^{\times n}$ with $\deg N_{[m]}(z) \leq m(\deg H-1)+1$.

$(c)$ If $H(z)$ is homogeneous  of degree $d\geq 2$, then,
 $N_{[m]}(z)$ is homogeneous of degree
$(d-1)m+1$ for any $m\geq 1$.
\end{corol}

\pf Again, by Theorem \ref{T4.3} 
and Lemma \ref{L5.1}, 
we know that the sequence 
$\{ N_{[m]}(z)\in \kzz |m\geq 1\}$
satisfies the recurrent relations 
Eqs.\,(\ref{M1}) and (\ref{Mm}).
If char.\,$k=0$, the corollary can be easily
proved by the mathematical induction 
on $m\geq 1$
via the recurrent relation 
Eq.\,(\ref{Mm}).
But, if char.\,$k=p>0$, 
the induction breaks down 
when $m\equiv 1$ $(\text{mod } p)$. 
However, we
can fix this problem as follows. 
Suppose the corollary holds 
for all  $1\leq l \leq kp$ 
for some $k\geq 1$. We 
consider $N_{[m]}(z)$ with $m=kp+1$.  
By Eq.\,$(\ref{L3.1.1-e2})$, we have
$$
H=N_t(z-tH)=\sum_{l\geq 1} t^{l-1} N_{[l]}(z-tH). 
$$
Comparing the coefficients of $t^{m-1}$ in the equation above, 
we have 
\begin{align}\label{Alt-Recurrent}
N_{[m]}(z)=-\res_t \, 
\sum_{l=1}^{m-1} t^{l-m-1} N_{[l]}(z-tH).
\end{align}
Note that, for any $1\leq l\leq m$,  
$\res_t \, t^{l-m-1} N_{[l]}(z-tH)$  
as the coefficient of 
$t^{m-l}$ of $N_{[l]}(z-tH)$
is obtained by replacing  
$(m-l)$ copies $z_i$'s by $(-H_i)$'s
in all possible ways 
for each monomial of $N_{[l]}(z-tH)$.
With this observation, it is easy to see that 
our mathematical induction arguments still can
go through at $m=kp+1$.
\epfv

Note that, by Corollary \ref{C5.4}, $(a)$, the infinite sum 
$\sum_{m=1}^\infty  N_{[m]}(z)t_0^{m-1}$ makes sense for any
$t=t_0\in k$. 
In particular, when $t=1$, $G_{t=1}(z)$ gives us 
the formal inverse $G(z)$ of $F(z)$. 
Now we can summarize the results above 
to formulate the following recurrent inversion formula.

\begin{theo} \label{T5.5} {\bf (Recurrent Inversion Formula)} 

Let $k$ be any integral domain of any characteristic. Let $H(z)$, 
$N_t(z)$ and $\{ N_{[m]}(z) |m\geq 1\}$ 
fixed as before. Then  

$(a)$ If $char.\, k=0$, $\{N_{[m]}(z) | m\geq 1 \}$ 
are completely determined  by 
\BQ
N_{[1]}(z) &=& H(z), \label{N1}\\
 N_{[m]}(z) &=&  \fr 1{m-1} \sum_{\substack{k+l=m\\ k, l\geq 1}} 
 \left [ N_{[k]}\pz \right ] N_{[l]}(z)  \label{Nm}
\EQ
for any $m\geq 2$.

$(b)$
If $char.\,k=p>0$, the recurrent relations above 
still hold for any $m\geq 2$ and $m\not \equiv 1$ 
$(\text{mod } p)$. 
When $m=kp+1$ for some $k\geq 1$, 
$N_{[m]}(z)$ can be obtained by 
Eq.\,$(\ref{Alt-Recurrent})$.
\end{theo}

When char.\,$k=p>0$, 
the inverse maps $G(z)$ 
can also be obtained by 
the following symbolic calculation.

\vskip2mm

\begin{algo}\label{Symbolic-Algo}
$(${\bf An Inversion Algorithm when char.\,${k=p>0}$}$)$
\vskip2mm
${\bf Step \,1:}$ 
Let $S$ be the set of the ordered triples $(i; I, J)$ 
with $1\leq i \leq n$ and 
$I, J \in (\bN^+)^{\times m}$ for some $m\geq 1$ 
such that the monomial 
$z_{i_1}^{j_1}z_{i_2}^{j_2}\cdots z_{i_m}^{j_m}$ appears
in $H_i(z)$
with a nonzero coefficient, say, $a_I^J(i)\in k$. 
Now let $A:=\{A_I^J(i) |(i; I, J) \in S\}$ 
be a set of free commutative variables 
and define $\wtilde F(z)\in \bZ [A]\langle\langle z\rangle \rangle^{\times n}$ 
 by replacing $a_I^J(i)$ 
by $A_I^J(i)$ in $F(z)$ for each triple 
$(i; I, J) \in S$. 

\vskip2mm

${\bf Step \, 2:}$ 
We view  $\wtilde F(z)$ as an automorphism of 
$\bZ [A]\langle\langle z\rangle \rangle$ 
over the base ring $\bZ [A]$
which is of characteristic zero.
Now we can apply 
the recurrent formulas Eqs.\,$(\ref{N1})$ 
and $(\ref{Nm})$ to calculate 
the inverse map 
$\wtilde G(z)$ of $\wtilde F(z)$.
Note that 
coefficients of all monomials  of $\wtilde G(z)$
are also in the base ring $\bZ[A]$.

\vskip2mm

${\bf Step\, 3:}$ 
To recover the inverse map $G(z)$ 
from $\wtilde G(z)$, we simply change all coefficients of 
$\wtilde G(z)$ 
by replacing each $A_I^J(i)$  
by $a_I^J(i)$ and 
each integer by its 
congruence class modulo $p$. 
\end{algo}

\begin{rmk}
In Step $1$ of the algorithm above, we may 
lift those coefficients $a_I^J(i)\in k$ 
which lie in $\bZ_p \subseteq k$ to any their pre-images 
in $\bZ$ instead to the corresponding 
formal variables $A_I^J(i)$. This may reduce the 
number of formal variables $A_I^J(i)$ involved 
and simplifies the algorithm substantially 
under certain circumstances.
\end{rmk}

Next let us consider the following example for the recurrent 
formula Eq.\,(\ref{Nm}) in Theorem \ref{T5.5}.

\begin{exam}\label{Nex-A}
Let $k$ be any integral domain with $char.\,k=0$
and $x, y$ two noncommutative free variables. 
Let $ad_y : k\langle\langle x, y\rangle\rangle\to 
k\langle\langle x, y\rangle\rangle$ be the $k$-linear 
map with $ad_y(u)=yu-uy$ for any 
$u\in k\langle\langle x, y\rangle\rangle$. 

Let $F(x, y)=(F_1(x, y), F_2(x, y))$ be 
the automorphism of the formal power series algebra 
$k\langle\langle x, y\rangle\rangle$ 
with  
\begin{align*}
\begin{cases}
F_1(x, y)&=x-(yx-xy)=(1-ad_y)(x),\\
F_2(x, y)&=y,
\end{cases}
\end{align*}
where $1$ 
in the expression denotes 
the identity map of 
$k\langle\langle x, y\rangle\rangle$. 

Now let us apply the recurrent formula in Theorem $\ref{T5.5}$
to determine the inverse map $G(x, y)$ of $F(x, y)$. 

Let $t$ be a central parameter and $F_t(x, y)$ be the 
special deformation discussed in Section $\ref{S4}$, i.e 
\begin{align*}
\begin{cases}
F_{t, 1}(x, y)&=x-t(yx-xy)=(1-t\, ad_y)(x),\\
F_{t, 2}(x, y)&=y.
\end{cases}
\end{align*}

Let $G_t(x, y)=z+t\,N_t(x, y)$ be the inverse map of $F_t(x, y)$. 
From the equation $F_{t, 2}(G_{t, 1}, G_{t, 2})=y$,
we see that $G_{t, 2}(x, y)=y$ and the second component 
of $N_t(x, y)$ must be zero. Therefore,
there exists $u_t(x, y)=\sum_{m\geq 1} t^{m-1} \, u_m(x, y) 
\in k\langle\langle x, y\rangle\rangle$ such that
\begin{align}
N_t(x, y)&=(u_t(x, y), 0), \\
N_{[m]}(x, y)&=(u_m (x, y), 0) 
\end{align}
for any $m\geq 1$.

With $N_{[m]}(x, y)$ having the form above, one can easily
 check that Eqs.\,$(\ref{N1})$ and $(\ref{Nm})$ 
in Theorem $\ref{T5.5}$ become
\begin{align}
u_{1}(x, y) &= xy-yx=ad_y(x), \label{u1} \\
u_{m}(x, y) &= \frac 1{m-1}  \sum_{\substack{k+l=m\\ k, l\geq 1}} 
 \left [ u_k(x, y) \frac{\p}{\p x} \right ] u_l (x, y). \label{um}
\end{align}
for any $m\geq 2$.
 
We claim that, for any $m\geq 1$, $u_m(x, y)=ad_y^m (x)$. This can 
be easily checked inductively by using Eqs.\,$(\ref{u1})$ 
and $(\ref{um})$ above along with the following 
simple observation: for any $k, l\geq 0$, we have
\begin{align}\label{NexA-E1}
\left [ad_y^k(x) \frac{\p}{\p x} \right ] ad_y^l(x) 
&=ad_y^l \left [ ad_y^k(x) \frac{\p}{\p x} \right ] x  \\
&=ad_y^l ad_y^k (x) \nno \\
&=ad_y^{k+l}(x). \nno 
\end{align}

Therefore, we have
\begin{align*}
G_{t, 1}(x, y)&=x+t N_{t, 1}(x, y)\\
&=x+t\sum_{m\geq 1} t^{m-1} \, ad_y^m(x)\\
&=\sum_{m\geq 0} t^{m} \, ad_y^m (x) \\
&=(1-t\,ad_y )^{-1} (x).
\end{align*}

In particular, the inverse map $G(x, y)$ of $F(x, y)$ is given by
\begin{align*}
\begin{cases}
G_1(x, y)&=(1-ad_y )^{-1} (x)=\sum_{m\geq 0}  ad_y^m (x),\\
G_2(x, y)&=y.
\end{cases}
\end{align*}

Note that the formula of $G_t(x, y)$ derived above 
can also be checked directly as follows.

\vskip2mm

\underline{\it Second Proof:} 
It is enough to check directly that 
\begin{align}
F_{t, 1}(G_{t, 1}, G_{t, 2})&=x, \\
F_{t, 2}(G_{t, 1}, G_{t, 2})&=y. \label{Ex-Eq1}
\end{align}

The second equation is obvious. Now consider the first one:
\begin{align*}
F_{t, 1}(G_{t, 1}, G_{t, 2})&=(1-t\, ad_y)\left (G_{t, 1}(x, y)\right ), \\
&=(1-t\, ad_y) \left ( (1-t\, ad_y)^{-1} (x) \right ), \\
&=x.
\end{align*}
\epfv
\end{exam}

\begin{rmk}
Considering the $($commutative$)$ Jacobian conjecture 
$($see \cite{BCW} and \cite{E}$)$, a naive noncommutative generalization 
of the Jacobian conjecture
would be: for any polynomial map $F$ of 
$k\langle\langle z\rangle\rangle$ with the Jacobian matrix 
$J F$ $($defined before Corollary $\ref{C2.2}$$)$
multiplicatively invertible, i.e. 
$JF\in GL_n(k\langle\langle z\rangle\rangle)$, it must be 
a polynomial automorphism of $k\langle\langle z\rangle\rangle$. 
The simple example above shows that this noncommutative 
generalization of the Jacobian conjecture is simply false, 
for the Jacobian matrix $J F_t(x, y)$ is
the $2\times 2$ identity matrix; while the inverse map $G_t(x, y)$ 
is not a polynomial map. For a correct noncommutative generalization
of the Jacobian conjecture, see \cite{Sc} and \cite{MSY}.
\end{rmk}

Finally, let us consider the following example 
for the symbolic Algorithm \ref{Symbolic-Algo} 
when the base ring $k$ has $char.\,k=p>0$.

\begin{exam}
Let $k$ be an integral domain of $char.\,k=p>0$ 
and $F(x, y)=(F_1(x, y), F_2(x, y))$ the polynomial map
of $k\langle \langle x, y\rangle \rangle$ given by
\begin{align*}
\begin{cases}
F_1 (x, y)&=x- s_0 (yx-xy)=(1-s_0\, ad_y)(x),\\
F_2 (x, y)&=y,
\end{cases}
\end{align*}
where $s_0$ is any fixed element of $k$.

Let $s$ be a formal variable. 
Applying Step $1$ of Algorithm $\ref{Symbolic-Algo}$, 
we get the polynomial map $\widetilde F_s (x, y)$
of  $\bZ[s]\langle \langle x, y \rangle \rangle$ 
with
\begin{align*}
\begin{cases}
\widetilde F_{s, 1}(x, y)&=x- s (yx-xy)=(1-s\, ad_y)(x),\\
\widetilde F_{s, 2}(x, y)&=y.
\end{cases}
\end{align*}

Now we consider Step $2$ of Algorithm $\ref{Symbolic-Algo}$. 
By the arguments in Example $\ref{Nex-A}$ above 
$($with the central parameter $t$ replaced by $s$$)$, 
we see that the inverse map $\widetilde G_s(x, y)$ of $\widetilde F_s(x, y)$ 
is given by 
\begin{align*}
\begin{cases}
\widetilde G_{s, 1}(x, y)&=\sum_{m\geq 0} s^m \, ad_y^m(x) ,\\
\widetilde G_{s, 2}(x, y)&=y.
\end{cases}
\end{align*}

Finally apply Step $3$, we get the inverse map $G(x, y)$ 
of $F(x, y)$ with 
\begin{align*}
\begin{cases}
G_1(x, y)&=\sum_{m\geq 0} s^m_0 \, ad_y^m(x) ,\\
G_2(x, y)&=y.
\end{cases}
\end{align*}

For example, when $k=\bZ_5=\bZ/5\bZ$ and $s_0=4=-1$, 
the inverse map $G(x, y)$ is given by
\begin{align*}
\begin{cases}
G_1(x, y)&=\sum_{m\geq 0} 4^m \, ad_y^m(x)=\sum_{m\geq 0} (-1)^m \, ad_y^m(x) ,\\
G_2(x, y)&=y.
\end{cases}
\end{align*}
\end{exam}

\renewcommand{\theequation}{\thesection.\arabic{equation}}
\renewcommand{\therema}{\thesection.\arabic{rema}}
\setcounter{equation}{0}
\setcounter{rema}{0}

\section{\bf An Expansion Inversion Formula by the
Planar Binary Rooted Trees}\label{S6}

In this section, we always assume the base ring $k$ 
has $char.\, k=0$.
We derive 
an expansion inversion formula
by the planar binary rooted trees for 
the inverse map $G(z)$ 
of automorphisms $F(z)$ of $\kzz$ 
(see Theorem \ref{T6.2}).  
Note that, unlike the tree expansion formula 
in \cite{BurgersEq}, which only holds for the
symmetric maps in commutative variables, 
the tree expansion inversion formula derived 
here works for all formal automorphisms 
in commutative or noncommutative variables.

First let us fix the following notations and conventions.

By a {\it rooted tree} we mean a finite
1-connected graph with one vertex designated as 
its {\it root}.
In a rooted tree
there are natural ancestral relations between vertices.  We say a
vertex $w$ is a child of vertex $v$ if the two are connected by an
edge and $w$ lies further from the root than $v$. 
We define the {\it degree}
of a vertex $v$ of $T$ to be the number of its children. 
A vertex is called a {\it leaf}\/ if it has no
children.  A rooted tree $T$ 
is said to be a {\it binary} 
if every non-leaf vertex of $T$ 
has exactly two children. 
A rooted tree $T$ 
is said to be a {\it planar} 
if the set of all children of 
each non-leaf vertex of $T$ is 
given a fixed linear order. 
A {\it planar rooted forest} is an ordered 
disjoint union of finitely 
many planar rooted trees.
A {\it planar binary rooted tree} is a rooted tree which is
both planar and binary.
When we speak of isomorphisms between rooted trees, 
we will always mean root-preserving isomorphisms. 

{\bf Notation:}

Once and for all, we fix the following notation for the rest of this paper.
\begin{enumerate}
\item We let $\mathbb T$ (resp.\,\,$\mathbb B$)
be the set isomorphism classes of all
rooted trees (resp.\,\,binary rooted trees). 
We denote by $\mathbb T^\cP$ (resp.\,\,$\mathbb B^\cP$)
the set of all planar rooted trees
(resp.\,\,planar binary rooted trees).
For any $m\ge1$, we let $\bT_m$, $\mathbb B_m$, $\mathbb T^\cP_m$ 
and $\mathbb B^\cP_m$
be the set of elements of  $\mathbb T$, $\mathbb B$,
$\mathbb T^\cP$ and $\mathbb B^\cP$, respectively, with $m$ vertices.

\item  We call  the rooted tree with one vertex the {\it singleton}, denoted by 
$\circ$.  For convenience, we also view the empty set
 as a rooted tree, denoted by $\emptyset$.

\item For any rooted tree $T$, we set the following notation:
\begin{itemize}
\item $\text{rt}_T$ denotes the root vertex of $T$.
\item $|T|$  denotes the
 number of the vertices of $T$ and 
$l(T)$ the number of leaves.
\item  $\widehat T$ denotes the rooted tree
obtained by deleting all the leaves of $T$. 
\end{itemize}
\end{enumerate}

For any set of rooted trees $T_1, T_2, ..., T_d$, 
we define $B_+(T_1, T_2,..., T_d)$ to 
be the rooted tree obtained by connecting all roots of  $T_i$ 
$(i=1, 2, ..., d)$ to a single new vertex, which is set to the root of  
the new rooted tree $B_+(T_1, T_2,..., T_d)$. For any rooted forest, say 
$T_1, T_2, ..., T_d$ ordered by their indices, we define 
$B_+(T_1, T_2,..., T_d)$ similarly, except we also order 
the set of children of the new root, which  
is set of roots of $T_i$'s, as the same order 
of $T_i$'s. Note that, for any 
$T_1, T_2 \in {\mathbb B}$, we have 
$B_+(T_1, T_2)\in {\mathbb B}$.

Next let us recall $T$-factorial $T!$ 
of rooted trees $T$, which was first introduced 
by D. Kreimer \cite{Kr}. It is defined inductively as follows. 
\begin{enumerate}
\item For the empty rooted tree $\emptyset$ and the singleton $\circ$, 
we set $\emptyset !=1$ and $\circ !=1$.
\item For any rooted tree $T=B_+(T_1, T_2, ..., T_d)$, 
we set 
\begin{align}\label{EE5.1}
T!=|T|\,T_1!\,T_2!\cdots T_d!.
\end{align}
\end{enumerate}

Note that, for the chains $C_m$ $(m\in \bN )$, i.e. 
the rooted trees with $m$ vertices and height $m-1$, we have $C_m!=m!$. 
Therefore the $T$-factorial of rooted trees can be viewed as a 
generalization of the usual factorial of natural numbers.

\begin{lemma}\label{L6.1}
$(a)$ For any non-empty binary rooted tree $T$, we have 
\begin{align}
|T|&= 2l(T)-1, \label{EE5.2}\\
|\widehat T|&=l(T)-1.\label{EE5.3}
\end{align}
$(b)$ For any $T \in \mathbb B^\cP$ with 
$T=B_+(T_1, T_2)$, we have
\begin{align}\label{E5.2}
\widehat T ! = (\ell(T)-1) \widehat T_1! \widehat T_2!
\end{align}
\end{lemma}
\pf $(a)$ can be 
proved easily 
by induction on the number of vertices. 
See Lemma $5.1$ 
in \cite{BurgersEq}, for example.
 
$(b)$ Note that, by the definition of the operation $B_+$
$\widehat T$, we have
$\widehat T=B_+(\widehat T_1, \widehat T_2)$. 
By Eqs.\,(\ref{EE5.1}) and (\ref{EE5.3}), we also have 
\begin{align}\label{E5.4}
 \widehat T!=|\widehat T|\, \widehat T_1!\, \widehat T_2!=
(l(T)-1)\, \widehat T_1!\, \widehat T_2!.
\end{align}
Hence we have Eq.\,(\ref{E5.2}).
\epfv

Now we fix an automorphism $F(z)=z-H(z)$ of $\kzz$ 
with $o(H(z))\geq 2$. Let $F_t(z)=z-tH(z)$ 
and $G_t(z)=z+tN(z)$ as in Section \ref{S4}.

We assign a $n$-sequence 
$N_T(z)\in \kzz^{\times n}$ 
for each non-empty planar binary 
rooted tree $T$ 
as follows.

\begin{enumerate}
\item For $T=\emptyset$, we set $N_T(z)=z$.
\item For $T=\circ$, we set $N_T(z)=H(z)$.
\item For any planar 
binary rooted tree $T=B_+(T_1, T_2)$, 
we set 
\begin{align}\label{B+N}
N_T(z)= \left [N_{T_1}(z)\pz\right] N_{T_2}(z).
\end{align}
\end{enumerate}

Now we are ready to state and prove 
the main theorem of this section.

\begin{theo} \label{T6.2}
For any $m\geq 1$, we have
\begin{align} \label{MainEq5.2}
N_{[m]}(z)
=\sum_{\substack {T\in \mathbb B^\cP_{2m-1}}} \frac 1{\widehat T!}  N_T(z)
=\sum_{\substack {T\in \mathbb B^\cP \\ l(T)=m} } \frac 1{\widehat T!}  
N_T(z).
\end{align}
Therefore, by Eq.\,$(\ref{Def-Nm})$ we have
\begin{align}
N_t(z)&=\sum_{T\in \mathbb B^\cP \backslash \emptyset} \frac {t^{l(T)-1}}{\widehat T!}  N_T(z),\\
G_t(z) & = \sum_{T\in \mathbb B^\cP} \frac {t^{l(T)}}{\widehat T!}  N_T(z).  
\end{align}
\end{theo}

\pf Note that, by Eq.\,(\ref{EE5.2}) in Lemma \ref{L6.1}, we have
\BQn
\mathbb B^\cP_{2m-1}&=&\{ T\in \mathbb B^\cP \, | \, l(T)=m\} \\
\mathbb B^\cP_{2m}&=& \emptyset,
\EQn 
for any $m\geq 1$. 
Hence the two sums in Eq.\,(\ref{MainEq5.2}) are equal to each other.

To prove Eq.\,(\ref{MainEq5.2}), we first set, for any $m\geq 1$, 
\begin{align*}
V_{[m]}(z)=\sum_{\substack {T\in \mathbb B^\cP_{2m-1}}} \frac 1{\widehat T!}  N_T(z).
\end{align*}
Then, by Theorem \ref{T5.5},
to show that $V_{[m]}(z)=N_{[m]}(z)$ for any $m\geq 1$, 
it will be enough to show that the sequence 
$\{V_{[m]}(z)\in \kzz  | m\geq 1\}$ also satisfies 
Eqs.\,(\ref{N1}) and (\ref{Nm}).

For the case $m=1$, since there is only one planar binary rooted tree $T$ 
with $l(T)=1$, namely, $T=\circ$, we have 
$V_{[1]}(z)=N_{T=\circ}(z)=H(z)$. 
Hence  Eq.\,(\ref{N1}) 
is satisfied.

For any $m\geq 2$, we consider
\allowdisplaybreaks{
\begin{align*}
&{} \quad \frac 1{m-1} 
\sum_{\substack {k, l\geq 1 \\ k+l=m}}
\left[ V_{[k]}(z)\pz\right] V_{[l]}(z) \\ 
&=
\sum_{\substack {T_1, T_2\in \mathbb B^\cP,\\ l(T_1)=k, l(T_2)=l,\\
k, l \geq 1,  k+l=m} } \frac 1{(m-1)\widehat T_1! \widehat T_2!}  \left[ N_{T_1}(z)\pz\right]  N_{T_2}(z) \\
&=
\sum_{\substack {T_1, T_2\in \mathbb B^\cP,\\ l(T_1)=k, l(T_2)=l,\\
k, l \geq 1,  k+l=m} } \frac 1{(m-1)\widehat T_1! \widehat T_2!}   N_{B_+(T_1, T_2)}(z) \\
\intertext{
Applying Eq.\,(\ref{E5.2}) in  Lemma \ref{L6.1}:}
&=
\sum_{\substack {T\in \mathbb B^\cP \\ l(T)=m} } \frac 1{\widehat T!}  
N_T(z)\\
&= V_{[m]}(z).
\end{align*} }
Hence we have Eq.\,(\ref{Nm}) 
for $V_{[m]}(z)$'s.
\epfv

Next let us consider the tree
expansion formula Eq.\,(\ref{MainEq5.2}) for the
polynomial map in Example \ref{Nex-A} 
(see Example \ref{Nex-C} below). As a bi-product,
we will get a proof for the following identities of 
the $T$-factorials $\widehat T!$ $(T\in \mathbb B^\cP)$.

\begin{propo}\label{P6.3}
For any $m\geq 1$, we have 
\begin{align}\label{T!-identities} 
\sum_{\substack {T\in \mathbb B^\cP_{2m-1}}} \frac 1{\widehat T!} = 
\sum_{\substack {T\in \mathbb B^\cP \\ l(T)=m} } \frac 1{\widehat T!}=1. 
\end{align}
\end{propo}

Note that the first equation in Eq.\,(\ref{T!-identities})
simply follows from the identity Eq.\,(\ref{EE5.2}).

\begin{exam}\label{Nex-C}
Let $F(x, y)$ and all the related 
notation as in Example $\ref{Nex-A}$.
Note that in this case $H(x, y)=(ad_y(x), 0)$. 
From Eq.\,$(\ref{B+N})$, it is easy to 
see inductively that, for any non-empty planar 
binary rooted tree $T$, the second component 
of $N_T(x, y)$ is also $0$. So we may write
$N_T(x, y)=(u_T(x, y), 0)$ for some 
$u_T(x, y)\in k\langle\langle x, 
y \rangle\rangle$. With this notation fixed, 
it is easy to check that, for any $T\in \mathbb B^\cP$ 
with $T=B_+(T_1, T_2)$, Eqs.\,$(\ref{B+N})$ and 
$(\ref{MainEq5.2})$ become respectively
\begin{align}
u_T(x, y)&= \left [u_{T_1}(x, y)\frac{\p}{\p x} \right] u_{T_2}(x, y), \label{B+u} \\
u_m(x, y)
&=\sum_{\substack {T\in \mathbb B^\cP_{2m-1}}} \frac 1{\widehat T!} \, u_T(x, y)
= \sum_{\substack {T\in \mathbb B^\cP \\ l(T)=m} } \frac 1{\widehat T!}\, u_T(x, y). 
\label{uT}
\end{align}

\vskip2mm

\underline{Claim:} for any $T\in \mathbb B^\cP$ with $T\neq \emptyset$,
we have
\begin{align}\label{NexC-E1}
u_T(x, y)=ad_y^{l(T)}(x).
\end{align}

\vskip2mm

\underline{Proof of Claim:} We use the mathematical induction 
on the number $l(T)$ of leaves of $T\in \mathbb B^\cP$.

First, when $T=\circ$, by the definition of $N_T(x, y)$, 
we know that $u_T(x, y)$ is the first component of $H(x, y)$,  
which is $ad_y(x)$. Hence Eq.\,$(\ref{NexC-E1})$ holds 
in this case.

Now assume Eq.\,$(\ref{NexC-E1})$ holds for any 
$T\in \mathbb B^\cP$ with $l(T)\leq m-1$ for 
some $m\geq 2$, and consider the case for 
$T\in \mathbb B^\cP$ with $l(T)=m$. 
Write $T=B_+(T_1, T_2)$ with 
$T_i\in \mathbb B^\cP$ $(i=1, 2)$ and $l(T_i)<m$. 
By Eq.\,$(\ref{B+N})$ and the induction assumption, 
we have
\begin{align*}
u_T(x, y)&=\left [u_{T_1}(x, y)\frac{\p}{\p x} \right] u_{T_2}(x, y)\\
&=\left [ad_y^{l(T_1)}(x)\frac{\p}{\p x} \right] ad_y^{l(T_2)}(x)\\
&=ad_y^{l(T_2)}(x)\left[ ad_y^{l(T_1)}(x)\right ] \\
&=ad_y^{l(T_1)+l(T_2)}(x)\\
&=ad_y^{l(T)}(x).
\end{align*}
\epfv

Now, from Eqs.\,$(\ref{uT})$ and $(\ref{NexC-E1})$, 
we see that, for any $m\geq 1$,
the first component $u_m(x, y)$ of 
$N_{[m]}(x, y)$ is given by 
\begin{align} \label{ExC-14}
u_m(x, y)
=\sum_{\substack {T\in \mathbb B^\cP_{2m-1}}} \frac 1{\widehat T!}  ad_y^{l(T)}(x)
=( \sum_{\substack {T\in \mathbb B^\cP_{2m-1}}}  \frac 1{\widehat T!})  \,
ad_y^{m}(x).
\end{align}

But from Example $\ref{Nex-A}$, we also know that 
\begin{align} \label{ExC-15}
u_m(x, y)=ad_y^{m}(x),
\end{align}
for any $m\geq 1$.

Comparing Eqs.\,$(\ref{ExC-14})$ and $(\ref{ExC-15})$, 
we get the identity Eq.\,$(\ref{T!-identities})$. 
\end{exam}

Finally, to make our arguments more complete and also 
the tree expansion formula Eq.\,(\ref{MainEq5.2})  
more convincing, let us end this paper 
with the following direct proof of Proposition \ref{P6.3}.

\vskip2mm

\underline{\it 2nd Proof of Proposition \ref{P6.3}:}
Let $s$ be a formal variable and $a(s)$ the following generating function 
\begin{align} 
a(s)=\sum_{\substack {T\in \mathbb B^\cP}} \frac 1{\widehat T!} s^{l(T)-1} 
=\sum_{m\geq 1} (\sum_{\substack {T\in \mathbb B^\cP \\ l(T)=m} } \frac 1{\widehat T!}) \, s^{m-1}. 
\end{align}

Consider
\allowdisplaybreaks{
\begin{align*} 
a^2(s)&=(\sum_{T_1\in \mathbb B^\cP} \frac 1{\widehat T_1!} s^{l(T_1)-1})
(\sum_{T_2\in \mathbb B^\cP} \frac 1{\widehat T_2!} s^{l(T_2)-1})\\
&=\sum_{(T_1, T_2) \in \mathbb B^\cP\times \mathbb B^\cP} 
\frac 1{\widehat T_1!\widehat T_2!}\, s^{l(T_1)+l(T_2)-2}\\
\intertext{Re-indexing the terms in the sum above 
by $T:=B_+(T_1, T_2)$ and noting that any $T\in \mathbb B^\cP$ with $|T|\geq 2$ 
can appear once and only once as $B_+(T_1, T_2)$ for some $T_i\in \mathbb B^\cP$:}
&=\sum_{\substack {T \in \mathbb B^\cP \\T=B_+(T_1, T_2) }} 
\frac 1{\widehat T_1!\widehat T_2!} \, s^{l(T_1)+l(T_2)-2} \\
\intertext{Applying Eq.\,$(\ref{E5.2})$:}
&=\sum_{\substack {T \in \mathbb B^\cP \\|T|\geq 2}} 
\frac {l(T)-1}{\widehat T!}\, s^{l(T)-2}\\
\intertext{Noting that when $|T|=1$, we have $T=\circ$ and $l(T)=1$:}
&=\sum_{T \in \mathbb B^\cP} 
\frac {l(T)-1}{\widehat T!} \, s^{l(T)-2}\\
&=\frac{d}{ds}\, a(s).
\end{align*} }
Therefore, we see that $a(s)$ satisfies the equations
\begin{align*} 
\begin{cases}
\frac{d\, a(s)}{ds}&=a^2(s),\\
a(0)&=1. 
\end{cases}
\end{align*}

But it is easy to check that the only formal power series 
solution of the equations above is $(1-s)^{-1}$. Therefore, we have 
$a(s)=(1-s)^{-1}$. By comparing the coefficients of $s^{m-1}$ 
$(m\geq 1)$ of $a(s)$ and $(1-s)^{-1}$, we get 
Eq.\,$(\ref{T!-identities})$. 
\epfv



\end{document}